\documentclass[a4paper,10pt,twoside]{amsart}
\usepackage{amsthm}
\usepackage{amsmath}
\usepackage{amssymb}
\usepackage{amsfonts}
\usepackage{amscd}
\usepackage[english]{babel}
\usepackage{stmaryrd}
\usepackage{enumerate}

\newtheorem{thm}{Theorem}
\newtheorem{cor}{Corollary}
\newtheorem{prop}{Proposition}
\newtheorem{rem}{Remark}

\newtheorem{exam2}{Example}
\newtheorem{conj}{Conjecture}

\newtheorem{defi2}{Definition}
\newtheorem{lem}{Lemma}

\newtheorem*{thm2}{Theorem}
\newtheorem*{prop2}{Proposition}

\begin{document}

    \title{Nonsymmetric generic matrix equations}
    \author{Gerald BOURGEOIS}
    
    \date{Sept-27-2014}
    \address{G\'erald Bourgeois, GAATI, Universit\'e de la polyn\'esie fran\c caise, BP 6570, 98702 FAA'A, Tahiti, Polyn\'esie Fran\c caise.}
    \email{bourgeois.gerald@gmail.com}
        

\begin{abstract} 
Let $(A_i)_{0\leq i\leq k}$ be generic matrices over $\mathbb{Q}$, the field of rational numbers. Let $K=\mathbb{Q}(E)$, where $E$ denotes the entries of the $(A_i)_i$, and let $\overline{K}$ be the algebraic closure of $K$. We show that the generic unilateral equation $A_kX^k+\cdots+A_1X+A_0=0_n$ has $\binom{nk}{n}$ solutions $X\in\mathcal{M}_n(\overline{K})$. Solving the previous equation is equivalent to solving a polynomial of degree $kn$, with Galois group $S_{kn}$ over $K$. Let $(B_i)_{i\leq k}$ be fixed $n\times n$ matrices with entries in a field $L$. We show that, for a generic $C\in\mathcal{M}_n(L)$, a polynomial equation $g(B_1,\cdots,B_k,X)=C$ admits a finite fixed number of solutions and these solutions are simple. We study, when $n=2$, the generic non-unilateral equations $X^2+BXC+D=0_2$ and $X^2+BXB+C=0_2$. We consider the unilateral equation $X^k+C_{k-1}X^{k-1}+\cdots+C_1X+C_0=0_n$ when the $(C_i)_i$ are $n\times n$ generic commuting matrices ; we show that every solution $X\in\mathcal{M}_n(\overline{K})$ commutes with the $(C_i)_i$. When $n=2$, we prove that the generic equation  $A_1XA_2X+XA_3X+X^2A_4+A_5X+A_6=0_2$ admits $16$ isolated solutions in $\mathcal{M}_2(\overline{K})$, that is, according to the B\'ezout's theorem, the maximum for a quadratic $2\times 2$ matrix equation.
\end{abstract}

\maketitle

    \section{Introduction} 
    S. Gelfand wrote in 2004 (cf. \cite{12}) : ``The problem of solving quadratic equations for matrices seems to be absolutely classical. It would be natural that such a problem should have been at least formulated, or even solved, in the 19th century at the latest. Still, I asked many people about this problem, and they directed me to various sources, but nowhere could I find even a mention of this problem''.\\
  \indent  Let $n\in \mathbb{N}_{\geq 2}$. In the present paper, we deal with polynomial equations where the coefficients are generic $n\times n$ matrices and the unknown is a $n\times n$ matrix; the underlying field is assumed to have characteristic $0$. Questions about generic matrices are solved in \cite{13} and \cite{14} or about formal matrices in \cite{16}. More generally, C. Procesi described, in \cite{15}, properties of the algebra of generic matrices.\\
    Let $\mathbb{Q}$ be the field of rational numbers. If $M$ is a $n\times n$ matrix, then $\chi_M$ denotes its characteristic polynomial, $\sigma(M)$ its spectrum and $\mathrm{tr}(M)$ its trace.
    \begin{defi2} (cf. \cite{33}) Let $\{a_{r;i,j}\;|\;1\leq i,j\leq n,1\leq r\leq k\}$ be independent commuting indeterminates over $\mathbb{Q}$ ; in other words, the $(a_{r;i,j})_{rij}$ are elements of a transcendental extension of $\mathbb{Q}$ and they are mutually transcendental over $\mathbb{Q}$. Then, when $r\in\llbracket 1,k\rrbracket$, the $n\times n$ matrices $A_r=[a_{r;i,j}]$ are called \emph{generic matrices} (over $\mathbb{Q}$) ; in the sequel, such matrices are assumed to be fixed. We consider the quotient field $K=\mathbb{Q}((a_{1;i,j})_{i,j},\cdots, (a_{k;i,j})_{i,j})$ and its algebraic closure $\overline{K}$.
   Let $f$ be a non-zero polynomial over $K$ in $k+1$ non-commuting indeterminates. We consider the so-called \emph{generic} matrix equation:
    \begin{equation} \label{def} f(A_1,\cdots,A_k,X)=0_n\text{ in the unknown }X=[x_{i,j}]\in\mathcal{M}_n(\overline{K}).\end{equation}
  $i)$ Assume that the previous equation has a finite positive number of solutions. If the entries of each solution can be calculated by radicals over $K$, then we say that Eq (\ref{def}) is solvable, else  we say that Eq (\ref{def}) is non-solvable.\\
  $ii)$ (cf. \cite{20}) A solution $X_0$ of Eq (\ref{def}) is called (geometrically) isolated if there is a neighborhood of $X_0$ that contains no other solution of the equation. 
    \end{defi2}
		Let $k,n\geq 2$ and $(A_i)_{0\leq i\leq n}$ be $n\times n$ genric matrices ; put $K=\mathbb{Q}((A_i)_i)$. In Section 2, we consider the unilateral equation of degree $k$ in the unknown $X\in\mathcal{M}_n(\overline{K})$
\begin{equation} \label{unilat}    A_kX^k+\cdots+A_1X+A_0=0_n.  \end{equation}
    Moreover, we study the nonsymmetric algebraic Riccati equation in $X\in\mathcal{M}_n(\overline{K})$
    \begin{equation} \label{comp} XAX+B_1X+XB_2+C=0_n,  \end{equation}
    where $A,B_1,B_2,C$ are $n\times n$ generic matrices.
       We reduce the study of Eq (\ref{comp}) to the following one 
    \begin{equation}  \label{final}     X^2+BX+C=0_n,     \end{equation}
    where $B,C$ are $n\times n$ generic matrices, that is to the generic Eq(\ref{unilat}) of degree $2$. \\			
We show
\begin{thm} \label{ric}
 $i)$ Eq (\ref{unilat}), in the unknown $X$, has exactly $\binom{kn}{n}$ simple solutions that are in $\mathcal{M}_n(K((\lambda_i)_i))$, where the $(\lambda_i)_{i\leq kn}$ are the roots of a polynomial of degree $kn$ that has $S_{kn}$ as Galois group over $K$. Moreover $K((\lambda_i)_i)$ is the minimal extension $L$ of $K$ such that $\mathcal{M}_n(L)$ contains all solutions of Eq (\ref{unilat}).\\
 $ii)$  Eq (\ref{unilat}) is solvable if and only if $n=k=2$. 
\end{thm}
    The previous result is part of the mathematical folklore, but, to the best of our knowledge, it has not been published and even proved. As an example, J.J. Sylvester (1885, cf. \cite[p. 234,235 and 196,197]{6}), studied Eq (\ref{unilat}) ; using the eigenvalue problem, he wrote, without any proof, that, exceptions apart, Eq (\ref{unilat}) admits $\binom{kn}{n}$ solutions. Moreover we show the following new result: solving Eq (\ref{unilat}) is equivalent to find the roots of a polynomial of degree $kn$ that has $S_{kn}$ as Galois group over $K$; in particular, the generic Riccati equation has $\binom{2n}{n}$ solutions and its resolution has the same complexity that solving a polynomial of degree $2n$ that has $S_{2n}$ as Galois group. \\
    \indent In Section $3$, we give two instances of specializations in $\mathbb{Q}$ of Eq (\ref{comp}) such that they are solvable until $n=4$. In particular, we study in details the generic equation $X^2+BX=0_n$.\\ 
\indent	Let $L$ be an algebraically closed field with characteristic $0$.	In Section 4, we prove 
		 \begin{thm}  \label{isol} Let $A_1,\cdots,A_k\in\mathcal{M}_n(L)$  be any known matrices and $B\in\mathcal{M}_n(L)$ be a known generic matrix. Let $g$ be a non-zero polynomial over $L$ in $k+1$ non-commuting indeterminates such that each monomial of $g$ contains explicitly the variable $X$. We consider the equation
  \begin{equation}   \label{isola}  \phi(X)=g(A_1,\cdots,A_k,X)=B\text{ in the unknown }X\in\mathcal{M}_n(L). 
  \end{equation}  
 Let $\delta$ be the degree of $\phi$ (see \cite[Section II-6.3]{22}). Then Eq (\ref{isola}) has generically (with respect to $B$) exactly $\delta$ solutions and all are simple.
  \end{thm}      
  \indent  To obtain explicit solutions, we use in Section $5$ and $7$, for low dimensions, Maple and Magma softwares ; this method permits to calculate the minimal polynomials of the entries of the solutions of Eq (\ref{comp}) and also to study non-unilateral quadratic equations. We consider the system of algebraic equations, in the unknowns $(x_{i,j})_{i,j}$ and the indeterminates $(b_{i,j})_{i,j},(c_{i,j})_{i,j}$, associated to Eq (\ref{final}). If $n=2$, then we completely solve this system. Yet, when $n=3,4$, we must use specializations in $\mathbb{Q}$ of $B,C$.
	\begin{defi2} We consider an Eq (\ref{def}) that has a finite number $\nu$ of solutions. Assume that finding the entries of each solution of Eq (\ref{def}) is equivalent to find the roots of a polynomial whose Galois group over $K$ is $G$. Then we denote by the solvability complexity (SC) of Eq (\ref{def}), the couple $(\nu,G)$. 
  \end{defi2}
    When $n\leq 3$, we study the non-unilateral generic equation	
    \begin{equation}  \label{plex1}   X^2+B_1XB_2+C=0_n.   \end{equation}
		Note that we cannot connect Eq (\ref{plex1}) to the eigenvalue problem. Then we know so little about the resolution of this type of equation.
    We show that, if $n=2$, then its SC is $(8,S_8)$. Note that we prove that there are $8$ simple solutions, without using software, but using only algebraic geometry considerations.		
		Moreover, when $n>2$, this equation admits a finite number of solutions and is non-solvable. \\
		When $n=2$, we consider the following specialization of Eq (\ref{plex1})
		\begin{equation} \label{plex2} X^2+BXB+C=0_2, \end{equation}
		where $B,C$ are $2\times 2$ generic matrices. With the Gr\"obner basis theory, we obtain an explicit solution of Eq (\ref{plex2}) and we show that its SC is $(6,S_6)$. \\
  \indent In Section $6$, we study the commuting unilateral matrix equations in the unknown $X\in\mathcal{M}_n(\overline{K})$
   \begin{equation} \label{commutat1}   X^k+X^{k-1}B_{k-1}+\cdots+XB_1+B_0=0_n,  \end{equation}
   \begin{equation} \label{commutat2}   X^k+B_{k-1}X^{k-1}+\cdots+B_1X+B_0=0_n.  \end{equation}
	When the commuting coefficients are generic (the precise definition is given in Section 6), we obtain the following:
	 \begin{thm}   \label{aston}
   We consider Eq (\ref{commutat1}) or Eq (\ref{commutat2}) where the $(B_i)_i$ are known $n\times n$ generic commuting matrices. Then there are $k^n$ solutions and any solution $X\in\mathcal{M}_n(\overline{K})$ commutes with the $(B_i)_i$.
  \end{thm}
\noindent	This result works for commuting Riccati Eq (\ref{comp}) or commuting Eq (\ref{plex1}) and does not work for the commuting Eq (\ref{plex2}): $X^2+BXB+C=0_2$, where $B,C$ are $2\times 2$ generic commuting matrices.\\
 \indent  The following theorem, proved in Section 7, gives, for $n=2$, an instance of a generic quadratic matrix equation that admits $16$ isolated solutions in $\overline{K}^4$, that is, according to the B\'ezout's Theorem (cf. \cite[Section 5.3]{18}), the maximal number of isolated solutions of a quadratic equation in $2\times 2$ matrices.
 \begin{thm}    \label{maxsol}
  The $2\times 2$ matrices $A,B_1,B_2,C,D,F$ are generic.
  We consider the equation in the unknown $X\in\mathcal{M}_2(\overline{K})$
  \begin{equation}   \label{degmax}   AXB_1X+XB_2X+X^2C+DX+F=0_2.
   \end{equation}
  $i)$ Then the $SC$ of Eq (\ref{degmax}) is $(16,S_{16})$.\\
		$ii)$ Moreover, the equation 
		\begin{equation}  \label{homog}  AXB_1X+XB_2X+X^2C=0_2 \end{equation}
		admits the sole solution $X=0_2$ with multiplicity $16$.
  \end{thm}
 \indent  In Section $8$, again when $n=2$, we present the different forms of the set of solutions of a non-generic Riccati equation.
     
  \section{The nonsymmetric unilateral algebraic equation}
   We begin with the Riccati equation. In particular,we will see that the complete Eq (\ref{comp}) is not solvable for $n\geq 3$. 
  Note that the form of Eq (\ref{comp}) is invariant by translation of the unknown.
   Let $X=Y+U$ where $U$ is a constant matrix to be chosen. We obtain $YAY+(UA+B_1)Y+Y(AU+B_2)+UAU+B_1U+UB_2+C=0_2$. Since $A$ is generic, it is invertible and we put $U=-A^{-1}B_2$. Thus the study of Eq (\ref{comp}) is reduced to the study of the following : let $A,B,C$ be $n\times n$ generic matrices. We consider the equation in the unknown $X\in\mathcal{M}_n(\overline{K})$
  \begin{equation} \label{princ}  XAX+BX+C=0_n.  \end{equation}
   Since $A$ is invertible, Eq (\ref{princ}) is equivalent to 
  $$(AX)^2+(ABA^{-1})(AX)+AC=0_n.$$
  Finally the study of Eq (\ref{comp}) is reduced to the study of the unilateral Eq (\ref{final}) in the unknown $X=[x_{i,j}]\in\mathcal{M}_n(\overline{K})$, where $B,C$ are $n\times n$ generic matrices.
	\subsection{The unilateral equation}
	More generally, we go to study the general unilateral equation ; we need the following result
\begin{thm2} (specialization) \cite[Section 5.8]{2} Let $P\in\mathbb{Q}(T)[X]=\sum_ia_i(T)X^i$ be an irreducible polynomial over $\mathbb{Q}(T)$, the set of rational fractions in $T=(t_1,\cdots,t_k)$ over $\mathbb{Q}$. If $t\in\mathbb{Q}^k$ is not a pole of one of the $(a_i)$, then we can associate the specialization $P_t \in\mathbb{Q}[X]$. Then the galois group of $P_t$ over $\mathbb{Q}$ is a subgroup of the Galois group of $P$ over $\mathbb{Q}(T)$.
\end{thm2}
 \begin{lem}    \label{gener}   Let $\Pi=[\pi_{i,j}]$ be a $p\times p$ matrix where the $(\pi_{i,j})_{i,j}$ are commuting indeterminates and $K=\mathbb{Q}((\pi_{i,j})_{i,j})$. Then $\chi_{\Pi}$ is irreducible over $K$ and its Galois group over $K$ is $S_p$. 
 \end{lem}
 \begin{proof}
 Let $P$ be a polynomial of degree $p$ with coefficients in $\mathbb{Q}$ that has $S_p$ as Galois group. We specialize $\Pi$ into $\Pi_0$ so that $\Pi_0$ is the companion matrix of $P$. According to the previous theorem, $S_p$ is a subgroup of the Galois group of $\chi_{\Pi}$ over $K$ and we are done. 
\end{proof}
Now we consider Eq (\ref{unilat}) :   $\;A_kX^k+\cdots+A_1X+A_0=0_n\;$. Note that, since $A_k$ is invertible, we may assume $A_k=I_n$. In conjunction with the eigenvalue problem of degree $k$, we need the following
\begin{thm2}   \cite[Section VIII-5]{28} Let $L $ be an algebraically closed field with characteristic $0$ and $(A_i)_{0\leq i\leq k}\in\mathcal{M}_n(L)$ be known matrices.
Let $X\in\mathcal{M}_n(L)$ be a solution of Eq (\ref{unilat}) and 
$$\phi(\lambda)=\det(\lambda^kA_k+\cdots+\lambda A_1+A_0)\in L[\lambda].$$
Then $\phi(X)=0_n$.
\end{thm2}
We prove \textbf{Theorem \ref{ric}}.
\begin{proof} 
Part 1.  Let $u$ be an eigenvector of $X$ associated to the eigenvalue $\lambda$. Then $(\lambda^kA_k+\cdots+\lambda A_1+A_0)u=0$ and, consequently, $\phi(\lambda)=0$.
Note that, since $\phi(\lambda)=\det(A_k)\lambda^{nk}+\cdots+\det(A_0)$, $\phi$ is a polynomial of degree $nk$ in $\lambda$. We consider the following specialization $\phi_0$ of $\phi$ in $\mathbb{Q}$
$$A_k=I_n;\text{ for every }2\leq i\leq k-1,\; A_i=0_n;\;A_1=[{a_1}_{i,j}],\text{ where }{a_1}_{i,j}=0\text{ except }$$
$${a_1}_{1,n}=(-1)^n;\;A_0=[{a_0}_{i,j}],\text{ where }{a_0}_{i,j}=0\text{ except }{a_0}_{i+1,i}=1,{a_0}_{1,n}=(-1)^n.$$
We obtain 
$$\phi_0(\lambda)=\det(\lambda^kI+\lambda A_1+A_0)=\lambda^{nk}-\lambda-1.$$
According to \cite{27}, the Galois group of $\phi_0$ over $\mathbb{Q}$ is $S_{nk}$ ; according to the specialization theorem, $\phi$ is irreducible over $K$ and its Galois group over $K$ is $S_{kn}$. Let $(\lambda_i)_{i\leq nk}$ denote the (simple) roots of $\phi$ ; the $n$ eigenvalues of $X$ are to be chosen amongst the $(\lambda_i)_i$.\\
\indent Part 2. We show that, if $\phi(\lambda_i)=0$, then $\mathrm{dim}(\ker({\lambda_i}^kA_k+\cdots+\lambda_i A_1+A_0))=1$. Otherwise $\psi(\lambda_i)=\det(\tau({\lambda_i}^kA_k+\cdots+\lambda_i A_1+A_0))=0$ where $\tau(U)$ denotes the $(n-1)\times (n-1)$ submatrix of a $n\times n$ matrix $U$ created from the rows $2,\cdots,n$ and the columns $1,\cdots,n-1$. Consequently, the resultant $\mathrm{result}(\phi(x),\psi(x))$ is zero, that is a polynomial relation, with coefficients in $K$, linking independant indeterminates. We obtain a contradiction if there is a specialization $\phi_0$ of $\phi$ in $\mathbb{Q}$ such that $\mathrm{result}(\phi_0(x),\psi_0(x))\not=0$. Again, we consider the specialization above ; then $\phi_0(x)=x^{nk}-x-1$ and $\psi_0(x)=1$ and, clearly, their resultant is not $0$. Thus, to each eventual eigenvalue $\lambda_i$ is associated a unique eigenvector $u_i$.\\
\indent Part 3. We exhibit $\binom{kn}{n}$ solutions of Eq (\ref{unilat}) and we show that there are no other solutions. We choose $J=\{i_1<\cdots<i_n\}$, a subset of $\{1,\cdots,nk\}$. Let $D_J=\mathrm{diag}(\lambda_{i_1},\cdots,\lambda_{i_n})$, $P_J=[u_{i_1},\cdots,u_{i_n}]$ and $X_J=P_JD_J{P_J}^{-1}$. Note that $X_J$ does not depend on the order chosen for the elements of $J$. One has $(A_kX^k+\cdots+A_0)P=0$, that implies that $X_J$ is a solution. It remains to prove that the eigenvalues of $X$ are simple. According to the previous theorem, if $X$ is a solution of Eq (\ref{unilat}), then $\phi(X)=0$. Therefore, the eigenvalues of $X$ are simple and $X$ is diagonalizable.\\
\indent Part 4. Assume that the solutions $(X_J)_J$ of Eq (\ref{unilat}) are in $\mathcal{M}_n(L)$ where $L\supset K$ ; we show, for instance, that $\lambda_1\in L$. Let $X_1,X_2$ be the solutions associated to $\{1<\cdots<n\}$ and $\{2<\cdots<n+1\}$. One has $\det(xI-X_1)=(x-\lambda_1)\cdots(x-\lambda_n)\in L[x]$ and $\det(xI-X_2)=(x-\lambda_2)\cdots(x-\lambda_{n+1})\in L[x]$. Then their $gcd$, $(x-\lambda_2)\cdots(x-\lambda_n)$ is in $L[x]$, and clearly, $\lambda_1\in L$.\\
Finally $ii)$ is an easy consequence of the previous four parts.
\end{proof}
 \begin{cor}
  The SC of Eq (\ref{unilat}) is $(\binom{kn}{n},S_{kn})$.
  \end{cor}
  \begin{proof}
  It is a consequence of Theorem \ref{ric} $i)$.   
  \end{proof}
\subsection{The particular case of the Riccati equation}
 Another consequence is that the generic algebraic Riccati equation, in dimension $n$, admits $\binom{2n}{n}$ solutions and its resolution is equivalent to solve a polynomial of degree $2n$ and of Galois group $S_{2n}$.\\
In fact, there is a theoretical solution of the non-generic Riccati equation. We give an outline of it, because we will use it in Section 3.
 Following the idea of W. Roth, J. Bell, J. Potter and B. Anderson (cf. \cite{30,26}), we associate to the non-generic Eq (\ref{comp}) the $2n\times 2n$ pseudo-hamiltonian matrix 
 $$M=\begin{pmatrix}-B_2&-A\\C&B_1\end{pmatrix}.$$
  If $X$ is  a solution of Eq (\ref{comp}), then the graph of $X$ is $M$-invariant (cf. \cite{10}). This essential property works because the LHS of Eq (\ref{comp}) is a $linear$ function of the coefficients. According to \cite{5}, there is a one-to-one correspondence between the set of solutions of Eq (\ref{comp}) and the set of $n$-dimensional $M$-invariant subspaces $E$ that satisfy the condition
 \begin{equation} \label{comple}  E \text{ is complementary to }\{0\}\times\overline{K}^n.\end{equation} 
 Such a subspace $E$ is spanned by $n$ vectors $f_1,\cdots,f_n$. Put $(f_1,\cdots,f_n)=\begin{pmatrix}U\\V\end{pmatrix}$ where $U,V$ are $n\times n$ matrices. Then the subspace $[f_1,\cdots,f_n]$ satisfies Condition (\ref{comple}) if and only if $U$ is invertible. Finally, the solution associated to $E$ is $V{U}^{-1}$, that does not depend on the choice of the basis $f_1,\cdots,f_n$. Note that if we know the eigenvalues of $M$, then we can calculate its Jordan form and consequently, we can solve Eq (\ref{comp}).
We deduce that,  when the non-generic Eq (\ref{comp}) has a finite number $\nu$ of solutions, then $\nu\leq \binom{2n}{n}$.
     \begin{rem}
$i)$ If $M$ is non-derogatory, then $M$ admits a finite number of invariant subspaces and the associated Eq (\ref{comp}) has a finite number (eventually $0$) of solutions (cf. \cite[ch. 17.8]{3}).\\
$ii)$ In \cite{4}, the authors study the non-generic equation 
\begin{equation} \label{hig}  AX^2+BX+C=0_n \end{equation}
 (where $A$ is not necessarily ivertible) and, using a method adapted from the analysis of the Riccati equation, they  construct solvents of  Eq (\ref{hig}). In the generic case, Eq (\ref{hig}) is equivalent to $X^2+B'X+C'=0$, that is Eq (\ref{final}). The conclusion is very different if we consider the following non-unilateral equation: 
$$X^2A+BX+C=0_n.$$
 Its LHS linearly depends on the coefficients ; yet, in the generic case, the previous equation is equivalent to $X^2+BXB'+C'=0$, that is Eq (\ref{plex1}), equation whose LHS is not a linear function of the coefficients (cf. Section 4).
 \end{rem}
   
    \section{Two specializations} 
    $\bullet$  Putting $A=I_n,B_1=B_2=B$, we specialize Eq (\ref{comp}) into    
        \begin{equation} \label{simple} X^2+BX+XB+C=0_n.  \end{equation} 
        where $B,C$ are $n\times n$ generic matrices.
      \begin{prop}
  Eq (\ref{simple}) admits $2^n$ distinct solutions in the unknown $X\in\mathcal{M}_n(\overline{K})$. If $n\leq 4$, then Eq (\ref{simple}) is solvable.  
  \end{prop}
  \begin{proof} 
  Eq (\ref{simple}) is equivalent to $(X+B)^2=B^2-C$. Since $C$ is generic, $B^2-C$ is generic and diagonalizable with non-zero distinct eigenvalues $(\lambda_i)_i$. There is an invertible matrix $P$ such that $B^2-C=P\;\mathrm{diag}(\lambda_1,\cdots,\lambda_n)P^{-1}$. Since $X+B$ and $B^2-C$ commute, $X=-B+P\;\mathrm{diag}(\mu_1,\cdots,\mu_n)P^{-1}$ where $\mu_i^2=\lambda_i$. If $n\leq 4$,then we can explicitly calculate the $(\lambda_i)_i$ and Eq (\ref{simple}) is solvable.
  \end{proof}
   $\bullet$ The equation 
   \begin{equation} \label{decomp} (XU+V)(RX+S)=0_n, \end{equation}
    where $U,V,R,S$ are $n\times n$ generic matrices, is a specialization of Eq (\ref{comp}) with the condition $C{B_2}^{-1}A=B_1$. This equation is studied in \cite{7}.\\
    Put $Y=RX+S$. Since $R,U$ are generic matrices, they are invertible and $X=R^{-1}(Y-S)$. We deduce easily that $(UY)^2+(-US+URVU^{-1})(UY)=0$. The problem is reduced to solving the equation in $Z\in\mathcal{M}_n(\overline{K})$
    \begin{equation} \label{binome} Z^2+TZ=0_n  \end{equation}
    where $T$ is a generic matrix. Note that Eq (\ref{binome}) plays an important role when solving Eq (\ref{comp}) (cf. \cite{8}). Clearly there are exactly $2^n$ solutions that commute with $T$. Yet, we go to see that there are many other solutions.\\
    \textbf{Notation} $i)$ Let $d$ be the Hilbert dimension of the polynomial ideal generated by Eq (\ref{decomp}) or by the associated Eq (\ref{binome}) ; let $\mathcal{S}$ be the union of irreducible components of dimension $d$ of the set of solutions of Eq (\ref{decomp}-\ref{binome}). \\
    $ii)$ Let $G(r,n)$ be the set of $r$-dimensional subspaces of $K^n$ (Grassmannian) ; it is a homogeneous space of dimension $\delta_r=r(n-r)$.
    \begin{prop} \label{grass}
  $i)$  If $n$ is even ($n=2p$), then $d=\delta_p=\left\lfloor{\dfrac{n^2}{4}}\right\rfloor$ and $\mathcal{S}$ consists of $\binom{2p}{p}$ affine varieties, whose equations have coefficients in $\mathbb{Q}[(\lambda_i)_i]$.\\
    If $n$ is odd ($n=2p+1$), then $d=\delta_p=\left\lfloor{\dfrac{n^2}{4}}\right\rfloor$ and $\mathcal{S}$ consists of $2\binom{2p+1}{p}$ affine varieties, whose equations have coefficients in $\mathbb{Q}[(\lambda_i)_i]$.\\
   $ii)$ If $n\leq 4$, then Eq (\ref{binome}) is solvable.     
    \end{prop}
    \begin{proof}
  $i)$  It is sufficient to consider Eq (\ref{binome}). Since $T$ is generic, we may assume $T=\mathrm{diag}((\lambda_i)_i)$ where the $(\lambda_i)_i$ are indeterminates. Here $M=\begin{pmatrix}0_n&-I_n\\0_n&T\end{pmatrix}$. Let $(e_i)_{i\leq 2n}$ be the canonical basis. Since $T$ is generic, $M$ is diagonalizable and $\ker(M)=[e_1,\cdots,e_n]$ is the eigenspace associated to the zero eigenvalue. An eigenvector associated to the eigenvalue $\lambda_i$ is $v_i=\dfrac{-1}{\lambda_i} e_i+e_{n+i}$. To construct a $M$-invariant $n$-vector space, we choose $r\leq n$ and $n-r$ vectors $v_{i_1},\cdots,v_{i_{n-r}}$ that generate a subspace $F_{n-r}$. Any subspace $G_r$ of $[e_1,\cdots,e_n]$ of dimension $r$ does not intersect $F_{n-r}$ and $F_{n-r}\oplus G_r$ is convenient. For the sake of simplicity, we take the vectors $v_1,\cdots,v_{n-r}$. According to Section $2$, the associated solutions are $VU^{-1}=\begin{pmatrix}D&Y\\0_{r,n-r}&0_{r,r}\end{pmatrix}$ where $D=\mathrm{diag}(-\lambda_1,\cdots,-\lambda_{n-r})$ and $Y$ is an arbitrary $(n-r)\times r$ matrix. Each choice of $G_r$ gives birth to a solution and, when $F_{n-r}$ is fixed, the set of solutions consists of an affine variety of dimension $\delta_r$. Now $\sup_r(\delta_r)$ is $\delta_p$ if $n$ is even and is $\delta_p=\delta_{p+1}$ if $n$ is odd. Since there are $\binom{n}{r}$ choices for $F_{n-r}$, we deduce the required results.\\
   $ii)$ If $n\leq 4$, then we can explicitly calculate the $(\lambda_i)_i$.
     \end{proof} 
     \begin{rem} The previous equation is in the form $Z^2+BZ+C=0_n$ where $BC=CB$ and $B^2-4C$ admits exactly $4$ square roots. According to \cite [p. 499]{4}, the solutions would be given by the usual formula for the roots of a scalar quadratic, that is, there would be exactly $4$ solutions. In fact, there is an infinity of solutions and many of them do not commute with $B,C$.
     \end{rem} 
    \section{A more general equation} 
    In this section, we show a general result that we shall need in the sequel.
    We consider Eq (\ref{def}); it is equivalent to a system of $n^2$ algebraic equations in the $n^2$ unknowns $(x_{k,l})_{k,l}$
    $$f(X)=0\;\Longleftrightarrow\;\text{for every }i,j\leq n\;,\;f_{i,j}(X)=0. $$
				\begin{defi2}  A solution $X_0$ of Eq (\ref{def}) is called singular if 
		$$\det(\dfrac{\partial f}{\partial X}(A_1,\cdots,A_k,X_0))=0,$$ 		
				that is $X_0$ is a critical point of the function $f(A_1,\cdots,A_k,.)$.  
					\end{defi2}
		\begin{exam2}  \label{crit}  Let $f:X\in\mathcal{M}_2(\mathbb{C})\rightarrow X^2+\begin{pmatrix}2&0\\0&-1\end{pmatrix}X$. Then
		$$\dfrac{\partial f}{\partial X}(X_0)=X_0\otimes I_2+I_2\otimes X_0^T+\begin{pmatrix}2&0\\0&-1\end{pmatrix}\otimes I_2.$$
		$i)$ The matrix equation $f(X)=\begin{pmatrix}-1&1\\0&2\end{pmatrix}$ admits the isolated solution\\
		$X_0=\begin{pmatrix}-1&1/3\\0&2\end{pmatrix}$. Then 
		$$\dfrac{\partial f}{\partial X}(X_0)=\begin{pmatrix}0&0&1/3&0\\1/3&3&0&1/3\\0&0&0&0\\0&0&1/3&3\end{pmatrix},$$
		and $X_0$ is a singular isolated solution.	\\
			$ii)$ The equation $f(X)=0_2$ admits the non-isolated solutions $X_{\alpha}=\begin{pmatrix}-2&\alpha\\0&0\end{pmatrix}$ when $\alpha\in\mathbb{C}$. Then $\dfrac{\partial f}{\partial X}(X_{\alpha})=\begin{pmatrix}-2&0&\alpha&0\\\alpha&0&0&\alpha\\0&0&-3&\alpha\\0&0&0&-1\end{pmatrix}$ and $\{X_{\alpha}\;|\;\alpha\in\mathbb{C}\}$ is a line of singular solutions.
				\end{exam2}
     \begin{defi2} (cf. \cite[Section 3.3]{18}) Let $X_0$ be an isolated solution of Eq (\ref{def}). Let $\mathcal{O}_{X_0}$ be the local ring in $X_0$ and $((f_{i,j})_{i,j})_{X_0}$ be the ideal generated by the $(f_{i,j})_{i,j}$ in $\mathcal{O}_{X_0}$. The multiplicity of $X_0$ is 
  $$\text{the dimension of }\mathcal{O}_{X_0}/((f_{i,j})_{i,j})_{X_0},\text{ as a }\overline{K}\text{-vector space}.$$
  \end{defi2}
   \begin{exam2}
  $i)$  Consider the solution $(0,0)$ of the system in $\mathbb{C}^2$ : $\{x^2=0,y^3=0\}$ ; the quotient  
     $\mathbb{C}[x,y]_{(0,0)}/(x^2,y^3)_{(0,0)}$ admits the basis $\{1,x,y,y^2,xy,xy^2\}$ and has dimension $6$. For the system $\{y^2-x^5=0,x^2-y^5=0\}$, the multiplicity of the solution $(0,0)$ is $4$ because it is locally equivalent to $\{x^2=0,y^2=0\}$.\\     
     If we know a Gr\"obner basis of the ideal generated by the $(f_{i,j})_{i,j}$, then fortunately, we can calculate (Maple does it) the total number of solutions of Eq (\ref{def})
 with multiplicity. For our second system, this number is $25$, that is $4$ solutions in $(0,0)$ and the $21$ solutions of the subsystem $\{y^{21}=1,x=y^{13}\}$.\\
$ii)$ The isolated solution in Example \ref{crit}, $i)$, has multiplicity $3$.
     \end{exam2} 
		Of course, there is a link between singular solution and multiplicity of solution.
		\begin{prop2}   (cf. \cite[Section 3.2]{20})
		Let $X_0$ be a solution of Eq (\ref{def}). Then $X_0$ is singular if and only if $X_0$ is an isolated multiple solution or $X_0$ is in an algebraic subset of solutions of positive dimension.		
		\end{prop2}
		Let $L$ be an algebraically closed field with characteristic $0$. Let $R_n$ denote the function field of $L^n$ with coefficients in $L$. We will need the next result; $i)$ corresponds to Sard's Theorem when $L=\mathbb{C}$ and $i),ii)$ constitute the second Bertini's Theorem (cf. \cite[Section II-6.2]{22}).
		\begin{thm2}  Let $f:L^n\rightarrow L^n$ be a polynomial function, with $f(L^n)$ Zariski-dense in $L^n$ ; let $\mathcal{S}$ denote the set of critical points of $f$ and $\delta=[R_n:f^*(R_n)]$ be the degree of $f$ (\cite[Section II-6.3]{22}). The following holds\\
		$i)$ There is a proper Zariski closed subset $V$ of $L^n$ such that $f(\mathcal{S})\subset V$.\\
		$ii)$ If $y$ is in the dense Zariski open set $L^n\setminus V$, then the fiber $f^{-1}(y)$ is at most $\delta$ points.\\
		$iii)$ There is a dense Zariski open set $U\subset L^n\setminus V$ such that, for every $y\in U$, the fiber $f^{-1}(y)$ is exactly $\delta$ points (From $L^n\setminus V$, one removes the ramification points of $f$, cf. \cite[Section II-6.3]{22}).
		\end{thm2}
		\begin{defi2}  \label{gen}
		Let $\mathcal{S}$ be a polynomial system  of $n$ equations in $n$ unknowns over a field $K$ such that its coefficients are polynomial functions of $\tau$ parameters $(u_i)_{i\leq \tau}\in K^{\tau}$. $\mathcal{S}$ is said to \emph{have generically the property }$\mathcal{P}$ if there is a Zariski open set $U\subset K^{\tau}$ such that if $(u_i)_{i\leq \tau}$ is chosen in $U$, then the property $\mathcal{P}$ is fulfilled.
		\end{defi2} 
		We use the rule, according to which, if a polynomial function $L^n\rightarrow L^n$ has a non-dense image, then its degree is $0$. \\ 
		Now we prove \textbf{Theorem \ref{isol}}.
  \begin{proof}
	Let $Z=\mathrm{im}(\phi)$.\\
	$Case\;1.$ $Z$ is not dense. Then there is a dense Zariski open set $U$ such that $Z\cap U$ is void ; generically (with respect to $B$), Eq (\ref{isola}) has no solutions.\\
	$Case\;2.$ $Z$ is dense. We use the previous Proposition and Theorem. By Theorem $ii)$, generically (with respect to $B$), the fibers of $\phi$ are finite. By Theorem $iii)$, generically, the number of solutions of Eq (\ref{isola}) is the degree of $\phi$. According to Proposition and Theorem $i)$, for a generic $B$, all solutions are simple.
   \end{proof}
			\begin{rem} $i)$ Clearly, Theorem \ref{isol} works for any algebraic function from $L^n$ to $L^n$.\\
$ii)$	In particular, Theorem \ref{isol} remains valid if $g$ is a polynomial in the $(A_i)_i$ and $X,X^T$. The generic equation $g(X)=XAX^T=B$ is not interesting since the image of $g$ is not Zariski-dense ; indeed, there is a solution if and only if $A$ and $B$ are congruent if and only if there are invertible matrices $P,Q$ such that $PAQ=B,PA^TQ=B^T$. Yet, random experiments seem to show that the generic so-called Riccati-like equation $g(X)=XAX^T+BX=C$ has $6$ solutions when $n=2$ and $36$ solutions when $n=3$. In the same way, the generic equation $XAX^T+X^TBX=C$ seems to have $4$ solutions when $n=2$ and $48$ solutions when $n=3$.
			\end{rem}
  \section{Random specializations in $\mathbb{Q}$ and formal calculus}  
 \subsection{Gr\"obner and Hilbert}    In this section and also in Section $7$, we show the interest of the formal calculus in the study of the generic matrix equations. Moreover, the Hilbert's irreducibility theorem  is usefull to construct, in a probabilistic way, a generic matrix equation.\\
	More precisely, we can consider Eq (\ref{def}) as a system of $n^2$ algebraic equations in the unknowns $(x_{i,j})_{i,j}$ and we seek, using the package ``Groebner'' of Maple, a Gr\"obner basis of the ideal generated by our $n^2$ equations (cf. \cite{17}) ; indeed, if one does not declare as variable the entries of the matrices $(A_i)_i$, then they are seen by Maple as indeterminates. Assume that we obtain such a Gr\"obner basis -that can be done only if $n$ is small, because the associated algorithms have a great complexity-. Then there are two cases\\
	Case 1. The Gr\"obner basis is in the form
	$$\text{There is }(i,j)\text{ such that }p(x_{i,j})=0\text{ where }p\in K[x]\text{ has degree }\delta\text{ and for every }$$
	$$(k,l)\not=(i,j)\;,\;x_{k,l}=p_{k,l}(x_{i,j}),\text{ where }p_{k,l}\in K[x]\text{ has degree }<\delta.$$
	Thus, there are $\delta$ solutions.
	\begin{defi2}
	We say that a solution of Eq (\ref{def}) is \emph{effective} when we are in the case 1, that is, we obtain an univariate polynomial $p$ such that the complexity of solving Eq (\ref{def}) is the same than the complexity of solving the equation $p=0$.	
	\end{defi2}
	Case 2. The calculation of the solutions is not effective. Yet we know that there is a finite number $\nu$ of solutions. 
The knowledge of a Gr\"obner basis allows to calculate $\nu$ (cf. \cite{31}) ; curiously, Maple does not do that when the polynomial equations contain some indeterminates, but it is not difficult to write a procedure that does the job.\\
 Note that the calculation is much longer in the first case. Assume that we are in 	Case 1 and let $G$ be the Galois group of $p$ over $K$.
   Recall that the specialization in $\mathbb{Q}$ of the indeterminates results in a Galois group over $\mathbb{Q}$ which is a subgroup of the original Galois group over $K$. Yet, if the choice is random, then (except if we are unlucky) we obtain the total group $G$. To study this fact, we need the \textbf{Hilbert's irreducibility theorem}
	\begin{thm2} \cite[Section 5.9]{2},\cite[Section 3.4]{32}. Let $P\in\mathbb{Q}(T)[X]=\sum_ia_i(T)X^i$ be an irreducible monic polynomial over $\mathbb{Q}(T)$ and $G$ be its Galois group over $\mathbb{Q}(t)$. For $\delta>0$, let 
	$$N(\delta)=\{t\in[0,\delta]^k\cap \mathbb{Z}^k\;|\;t\text{ pole of some }a_i\text{ or the Galois group of  }P_t  \text{ is not }G\}.$$	
	Then the cardinality of $N(\delta)$ is in $O(\delta^{k-1/2+\epsilon})$ for every $\epsilon>0$.
	\end{thm2}	
   The previous result shows that a random specialization in $\mathbb{Q}^k$ of the indeterminate $T$ allows to retrieve the total Galois group of the polynomial $P$. 
  If $\delta$ is not a large number, then the probability of falling on a $t$ that is not convenient, is not negligible, as we see about the following test where the entries are randomly chosen in $\llbracket -2,2\rrbracket$.  
   \begin{rem}
  When $n=3$, there exist solvable instances $(B,C)$ of the non-generic Eq (\ref{final}), such that $BC-CB$ is not nilpotent (cf. Theorem \ref{aston} when $B,C$ commute), as this one
  $$B=\mathrm{diag}(1,-1,2)\;,\;C=\begin{pmatrix}2&-1&2\\-1&0&1\\1&2&0\end{pmatrix}.$$
  Here, $R=R_4R_6$, where $R_4,R_6$ are irreducible polynomials of degrees $4,6$ and $R_6$ is solvable.
  \end{rem}
	Unfortunately, there is no result that gives explicitly an integer $k(P)$ such that if $k$ random specializations $P_0$ in $\mathbb{Q}$ of $P$ all have $G$ as Galois group over $\mathbb{Q}$, then the probability that the Galois group of $P$ over $K$ is $G$, is (for instance) approximately $0.99$.
  \subsection{Tests on the generic Riccati equation}
	Let $n\in\llbracket 2,4\rrbracket$. According to Section $2$, it is sufficient to consider Eq (\ref{final}) ; we know its SC and, consequently, using the  Gr\"obner basis theory, we seek an effective solution. Since $B$ is generic, we may assume that $B=\mathrm{diag}((\lambda_i)_i)$, where the $(\lambda_i)_i$ are distinct ( the $(\lambda_i)_i$ can be calculated by radicals from the entries of $B$). Let $C=[c_{i,j}]$. Thus we must solve a system of $n^2$ algebraic equations in the $n^2$ unknowns $(x_{i,j})_{i,j}$ and in the $n^2+n$ indeterminates $(\lambda_i)_i,(c_{i,j})_{i,j}$. Let $L$ be the field $\mathbb{Q}((\lambda_i)_i,(c_{i,j})_{i,j})$ and $t$ be the trace of $X$. Since there are a finite number of solutions in $X$, then the Hilbert dimension of the ideal generated by these polynomials is equal to $n^2+n$, the number of indeterminates. The reason that we choose the unknown $t$ is the following : the minimal polynomial over $L$ of $t$ is in the form $P(z)=z^\tau+\sum_{i=0}^{\tau-1}a_iz^i$ where $\tau=\binom{2n}{n}$. At least when $n\leq 4$, $Q(u)=P(u-\dfrac{a_{\tau-1}}{\tau})$ is an even polynomial. Therefore $Q(u)=R(u^2)$ where $R$ is a polynomial of degree $\tau/2$ ; since to each value of $t$ is associated a solution of Eq (\ref{final}), the problem is reduced to solving a polynomial of degree $\frac{1}{2}\binom{2n}{n}$. The results are as follows\\
$\bullet$ The calculation is effective only when $n=2$. Here $R$ has degree $3$ and we deduce explicitly the values of $t$.\\
$\bullet$	When $n=3$, too much memory is required, even with a cluster. Thus, when $n=3,4$ we must specialize the matrices $B,C$. Then we obtain effective solutions. As expected, for almost all tests, the Galois group of the minimal polynomial of $t$ is $S_{2n}$.\\ 
\indent	In the sequel of the Section, we consider non-unilateral quadratic equations. There are few results concerning these equations. In particular the authors of \cite{25}, giving a general version of noncommutative Vieta theorem, consider only the unilateral matrix equations.
  \subsection{Tests on Eq (\ref{plex1})}
	Now we consider Eq (\ref{plex1}) : $X^2+BXC-D=0_n$, in the unknown $X\in\mathcal{M}_n(\overline{K})$,
    where $B,C,D$ are $n\times n$ generic matrices. Unlike to Eq (\ref{comp}), its LHS does not linearly depend on the coefficients and, consequently, we do not know any explicit form of the solutions of Eq (\ref{plex1}) ; thus we study this equation using a specific method when $n=2$ and formal calculus when $n=3$. We obtain that Eq (\ref{plex1}) has a larger SC than that of Eq (\ref{comp}). 
  \begin{prop}  \label{pato}
  When $n=2$, the SC of Eq (\ref{plex1}) is $(8,S_8)$. For any $n> 2$, Eq (\ref{plex1}) has a finite number of solutions and is non-solvable. 
  \end{prop}
  \begin{proof}  According to Theorem \ref{isol}, Eq (\ref{plex1}) has a finite number of solutions.\\
  \textbf{Part 1}. $n=2$. There are $8$ solutions.
	 The proof, due to D. Barlet, uses algebraic geometry arguments.\\
	$i)$ The matrices $B,C$ admit $n$ distinct eigenvalues $(\lambda_i)_i,(\mu_i)_i$. The tensor products 
	$$B\otimes I:Y\in\mathcal{M}_n(K)\rightarrow BY,I\otimes C^T:Y\in\mathcal{M}_n(K)\rightarrow YC$$
	commute, are diagonalizable and admit the eigenvalues  $(\lambda_i)_i$, $n$ times each, and $(\mu_i)_i$, $n$ times each. If $p\leq n$, then the endomorphism 
	$$\psi_p:w\in \mathcal{M}_{1,n}(K)\rightarrow w(C-\mu_p I_n)\text{ has rank }n-1;$$
	let $w^p=(w_i^p)_i$ be a non-zero vector of $\ker(\psi_p)$. We need the following
	\begin{lem}
	Let $B,C$ be generic matrices with spectra $(\lambda_i)_i,(\mu_i)_i$ and $(w^p)_p$ be the vectors previously associated to $C$. Then there is a basis $(E_{i,j})_{i,j}$ of $\mathcal{M}_n(K)$ such that
	$$\text{for every }i,j\;,\; BE_{i,j}=\lambda_iE_{i,j}\;,\; E_{i,j}C=\mu_jE_{i,j}\text{ and }E_{i,j}E_{h,k}=w_h^jE_{i,k}.$$	
	\end{lem}
	\begin{proof}
	We may assume $B=\mathrm{diag}((\lambda_i)_i)$. For every $i,j$, let $E_{i,j}$ be the matrix such that each row is $0$, except the $i^{th}$ row, which is $w^j$. Clearly, these $(E_{i,j})_{i,j}$ constitute a basis and we are done.
	\end{proof}
	Now $n=2$, $B=\mathrm{diag}(\lambda_1,\lambda_2)$ and, choosing the eigenvectors of $B$, we may assume $C=\begin{pmatrix}\mu_1+\nu_1&\mu_1\nu_1-\mu_2\nu_1+{\nu_1}^2	\\-1&\mu_2-\nu_1\end{pmatrix}$.\\
	Using the previous lemma, we obtain 
	$$E_{1,1}=\begin{pmatrix}1&\nu_1\\0&0\end{pmatrix},E_{1,2}=\begin{pmatrix}1&\nu_2\\0&0\end{pmatrix},E_{2,1}=\begin{pmatrix}0&0\\1&\nu_1\end{pmatrix},E_{2,2}=\begin{pmatrix}0&0\\1&\nu_2\end{pmatrix},$$
	where $\mu_1+\nu_1=\mu_2+\nu_2$. Put $X=\sum_{i,j}z_{i,j}E_{i,j}$ and $D=\sum_{i,j}d_{i,j}E_{i,j}$. Thus Eq (\ref{plex1}) is equivalent to the system of $6$ equations  $(\Sigma_i)_{i\leq 6}$ in the unknowns $(z_{i,j})_{i,j},\alpha,\beta$
	$$z_{1,1}(\alpha+\lambda_1\mu_1)-\nu_2\beta-d_{11}=0\;,\;z_{1,2}(\alpha+\lambda_1\mu_2)+\nu_1\beta-d_{1,2}=0,$$
	$$z_{2,1}(\alpha+\lambda_2\mu_1)+\beta-d_{2,1}=0\;,\;z_{2,2}(\alpha+\lambda_2\mu_2)-\beta-d_{2,2}=0,$$
	$$\alpha=z_{1,1}+z_{1,2}+\nu_1z_{2,1}+\nu_2z_{2,2}\;,\;\beta=z_{1,1}z_{2,2}-z_{1,2}z_{2,1}.$$
	Let $\phi:(\alpha,\beta)\in {\overline{K}}^2\rightarrow$ 
	$$(z_{1,1}=\dfrac{d_{1,1}+\nu_2\beta}{\alpha+\lambda_1\mu_1}\;,\;z_{1,2}=\dfrac{d_{1,2}-\nu_1\beta}{\alpha+\lambda_1\mu_2}\;,\;z_{2,1}=\dfrac{d_{2,1}-\beta}{\alpha+\lambda_2\mu_1}\;,\;
	z_{2,2}=\dfrac{d_{2,2}+\beta}{\alpha+\lambda_2\mu_2})\in {\overline{K}}^4$$ 
	and $S$ be the surface $\mathrm{im}(\phi)$.\\
$ii)$	The Zariski closure of $S$ in the projective space $\mathbb{P}_4$ has degree $4$. Indeed, from the equations $\Sigma_1,\Sigma_2,\Sigma_5,\Sigma_6$, we calculate $\alpha,\beta$ as functions of $z_{1,1},z_{1,2}$ ; then $\Sigma_3,\Sigma_4$ can be rewritten as equations of two generic quadrics $Q_1,Q_2$ in the coordinates $(z_{i,j})_{i,j}$. Since $S=Q_1\cap Q_2$, the degree of the Zariski closure of $S$ is the product of the degrees of $Q_1,Q_2$, that is $4$. \\	
$iii)$	The Zariski closure of the graph $\Gamma$ of $\phi$ has degree $4$ in $\mathbb{P}_6$. Indeed, the generic affine $2$-plane of $\overline{K}^4$ is the trace of the generic affine $4$-plane of $\overline{K}^6$, and the degree of $\Gamma$ is the cardinality of the intersection of $\Gamma$ with a generic $4$-plane ; clearly $\Gamma$ has dimension $2$.\\
 We consider the intersection, in $\overline{K}^6$, of the hyperplane $H$ and the quadric $Q$ defined as follows
	$$H\;:\;\alpha=z_{1,1}+z_{1,2}+\nu_1z_{2,1}+\nu_2z_{2,2}\;,\;Q\;:\;\beta=z_{1,1}z_{2,2}-z_{1,2}z_{2,1}.$$
	Clearly $H\cap Q$ has degree $2$ and dimension $4$. The set of solutions of our system is $\Gamma\cap (H\cap Q)$. The sum of dimensions of $H$ and $H\cap Q$ is $6$ ; according to B\'ezout's theorem, the number of solutions, in the generic case, is the product of the degrees, that is $8$, and these solutions are simple ; this is also the number of solutions of Eq (\ref{plex1}).\\
	\textbf{Part 2}. $n=2$. The SC is $(8,S_8)$.\\
	We may assume that $B$ is diagonal. With Maple, we obtain a Gr\"obner basis ; unfortunately it is non-effective, but one can deduce (cf. Section 5.1) that the number of solutions is $8$ and, thus, we find again the result of Part 1. Let $p\in K[x]$ be the monic polynomial, the roots of which, are the solutions in $x_{1,1}$ with multiplicity ; clearly $\mathrm{degree}(p)\leq 8$. A specialization of $B,C,D$ in $\mathbb{Q}$ gives an effective Gr\"obner basis directed by $p_0(x_{1,1})=0$, where $p_0$ has degree $8$ and $S_8$ as Galois group over $\mathbb{Q}$.	
	Since $p_0$ is irreducible, $\mathrm{degree}(p)=8$, $p_0$ is a specialization of $p$ in $\mathbb{Q}$ and $p$ is irreducible over $K$. If $G$ is the Galois group of $p$ over $K$, then $S_{8}$ is a subgroup of $G$ and $G=S_{8}$. Finally, $p$ has exactly $8$ simple roots and, to each root $\tilde{x}_{1,1}$, is associated exactly one solution $(\tilde{x}_{i,j})_{i,j}$ of  Eq (\ref{plex1}). That implies that the $(\tilde{x}_{i,j})_{(i,j)\not=(1,1)}$ are in $K[\tilde{x}_{1,1}]$ and solving  Eq (\ref{plex1}) is equivalent to solve $p=0$.  \\
 \textbf{Part 3}.  $n>2$.\\
 Eq (\ref{final}) is a specialization of Eq (\ref{plex1}). According to Theorem \ref{ric},  Eq (\ref{final}) has a finite number of solutions and is non-solvable ; thus Eq (\ref{plex1}) is non-solvable.
    \end{proof}
 	When $n=3$, almost all the specializations of $B,C,D$ in $\mathbb{Q}$ give an effective Gr\"obner basis directed by $p_0(x_{1,1})=0$, where the  polynomial $p_0$ has degree $56$ and, according to the software Magma, $S_{56}$ as Galois group over $\mathbb{Q}$. Thus, considering the Hilbert's theorem, we conjecture the following
  \begin{conj}
  When $n=3$, the SC of Eq (\ref{plex1}) is $(56,S_{56})$.
  \end{conj}
	 	Now we consider Eq (\ref{plex2}): $X^2+BXB+C=0_n$, a specialization of Eq (\ref{plex1}), where $B,C$ are generic matrices.
	\subsection{Tests on Eq (\ref{plex2})}
	$\bullet$ $n=2$. We may assume that $B=\mathrm{diag}(b_1,b_2)$ and we put $C=[c_{i,j}],X=[x_{i,j}]$. We consider an algebraic system of $4$ quadratic equations in $4$ unknowns and $6$ indeterminates. With Maple, we obtain a Gr\"obner basis in the form 
	$$p(x_{1,1})=0,\text{ where }\mathrm{degree}(p)=6\text{ and for every }$$
	$$(i,j)\not=(1,1)\;,\;x_{i,j}=p_{i,j}(x_{1,1})\text{ where }\mathrm{degree}(p_{i,j})\leq 5,$$
	that is, an effective solution of Eq (\ref{plex2}). A specialization of $B,C$ in $\mathbb{Q}$, gives a polynomial $p_0$ that has $S_6$ as Galois group over $\mathbb{Q}$. Then the Galois group of $p$ over $K$ is also $S_6$. \\	
	$\bullet$ $n>2$. Assume that there is $n>2$ such that Eq (\ref{plex2}) is solvable. We choose $B=\mathrm{diag}(B',0_{n-2}), C=\mathrm{diag}(C',0_{n-2})$ where $B',C'$ are $2\times 2$ generic matrices.  There are $6$ particular solutions of Eq (\ref{plex2}) in the form $X=\mathrm{diag}(X',0_{n-2})$ satisfying $X'^2
+B'X'B'+C'=0_2$.	These solutions cannot be written with radicals, that is contradictory. Thus we can conclude with the proposition
\begin{prop}
	When $n=2$, we obtain an effective solution of Eq (\ref{plex2}) and its SC is $(6,S_6)$. When $n>2$, Eq (\ref{plex2}) admits a finite number of solutions and is non-solvable.
	\end{prop}
  \section{When the coefficients commute}
  Now we assume that the coefficients in our equation commute. 
	\begin{defi2}
	Let $B_0=[b_{0;i,j}],\cdots,B_{k-1}=[b_{k-1;i,j}]\in\mathcal{M}_n(L)$ where $L$ is a transcendental extension of $\mathbb{Q}$ ; let $\mathcal{I}$ be the ideal over $\mathbb{Q}$ generated by the relations of commutativity between the $n\times n$ matrices $U_0,\cdots,U_{k-1}$: 
	$$\text{ for every }0\leq i< j\leq k-1, U_iU_j=U_jU_i.$$
	We assume that the $(B_r)_r$ are commuting matrices and that, for every polynomial $P\notin \mathcal{I}$ over $\mathbb{Q}$, $P(\{b_{r;i,j}\}_{rij})\not=0$. Then the $(B_r)_r$ are said \emph{generic commuting matrices} over $\mathbb{Q}$.	                                         
	\end{defi2}	
	Let $B_0,\cdots,B_{k-1}$ be $n\times n$ generic commuting matrices and let\\
	$K=\mathbb{Q}(\{b_{0;i,j}\},\cdots,\{b_{k-1;i,j}\})$. We consider the unilateral Eq (\ref{commutat1}) and Eq (\ref{commutat2}) in the unknown $X\in\mathcal{M}_n(\overline{K})$.
      Note that, in the previous equations, the LHS is a linear function of the generic coefficients and the monomials of type $BX^k$ and $X^lC$ are not mixed.\\
 We show \textbf{Theorem \ref{aston}}. 
  \begin{proof}	
   By transposition, solving Eq (\ref{commutat2}) is reduced to solving Eq (\ref{commutat1}). Thus we study only Eq (\ref{commutat1}). Remark that the discriminant of $\chi_{B_0}$ is not in $\mathcal{I}$; to see that, it suffices to consider a specialization of the $(B_r)_r$ in the form $B_0=\mathrm{diag}(u_1,\cdots,u_n)$ where the $(u_i)$ are distinct and for every $r\geq 1$, $B_r=I_n$. Thus $B_0$ is similar to a diagonal matrix $\mathrm{diag}(\lambda_1,\cdots,\lambda_n)$ whose distinct eigenvalues are conjugate over $K$ and the commutant of $B_0$ is $K[B_0]$. Thus, for every $1\leq j\leq k-1$, there exists a unique polynomial $P_j$ of degree $n-1$ and with coefficients in $K$ such that $B_j=P_j(B_0)$.	
	Note that the entries of $B_0$ and the coefficients of the $\{P_j\}_{1\leq j\leq k-1}$ constitute independent (over $\mathbb{Q}$) indeterminates in $K$.	
	Let $(e_i)_{i\leq n}$ be a basis consisting of eigenvectors of $B_0$. Thus for every $i\leq n$
  $$(X^k+P_{k-1}(\lambda_i)X^{k-1}+\cdots+P_1(\lambda_i)X+\lambda_iI_n)e_i=0.  $$
  Put 
  $$\theta(x,\lambda)=x^k+P_{k-1}(\lambda)x^{k-1}+\cdots+P_1(\lambda)x+\lambda\;,\;P(\lambda)=\{P_r(\lambda)\}_{1\leq r\leq k-1}$$
   and for every $i$, let $\mu_{i,1},\cdots,\mu_{i,k}$ be the roots in the unknown $x$ of $\theta(x,\lambda_i)$. \\
  If the $(\mu_{i,j})_{i\leq n,j\leq k}$ are not distinct, then there is $i\leq n$ such that the discriminant
	$$\mathrm{discrim}(\theta(x,\lambda_i),x)=\phi(P(\lambda_i)) \text{ is zero, }$$
  or there are $p<q\leq n$ such that the resultant
  $$ \mathrm{result}(\theta(x,\lambda_p),\theta(x,\lambda_q),x)=\psi(P(\lambda_p),P(\lambda_q)) \text{ is zero. }$$
  Note that $\phi,\psi$ are polynomials with coefficients in $\mathbb{Q}$.	
	We consider the symmetrizations of $\phi(P(\lambda_i))$ and $\psi(P(\lambda_p),P(\lambda_q))$:
	$$\Pi_{\sigma\in S_n}\phi(P(\lambda_{\sigma(i)}))=0\text{ and }\Pi_{\sigma\in S_n}\psi(P(\lambda_{\sigma(p)}),P(\lambda_{\sigma(q)}))=0,$$	
	that are polynomial relations, with coefficients in $\mathbb{Q}$, linking the entries of $B_0$ and the coefficients of the $(P_j)_{1\leq j\leq k-1}$. We obtain a contradiction if we show that $\phi(\lambda_i)$ and $\psi(\lambda_p,\lambda_q)$ are not identically zero.\\ 
 \indent We specialize the $(B_i)_i$ into the $(B'_i)_i\in{\mathcal{M}_n(\mathbb{Q})}^k$, putting, for every $r\leq n$, $\lambda_r=r$ and, for every $0\leq s\leq k-1$, $P_s(x)=x^{s+1}\; \mod\; \chi_{B_0}(x)$. Then a calculation shows that the associated $(\mu'_{i,j})_{i\leq n,j\leq k}$ are distinct. Note that every $k$-tuple of commuting matrices in a neighborhood of $(B'_i)_i$ in ${\mathcal{M}_n(\mathbb{Q})}^k$ satisfies the previous property.\\
 \indent  We return to the case where the $(B_i)_i$ are generic commuting matrices. There are distinct $\mu_{i,1},\cdots,\mu_{i,k}$ such that $(X-\mu_{i,1}I_n)\cdots(X-\mu_{i,k}I_n)e_i=0$. Thus there is $j_i\leq k$ such that $\mu_{i,j_i}\in\sigma(X)$. Since the $(\mu_{i,j})_{i\leq n,j\leq k}$ are distinct, $X$ has the following simple eigenvalues : $(\mu_{i,j_i})_{i\leq n}$. We deduce that, for every $i$, $Xe_i=\mu_{i,j_i}e_i$ and $X$ commute with the $(B_i)$. It remains to solve $n$ univariate polynomials of degree $k$. 
  \end{proof}
   \begin{rem} An essential argument in the previous proof is that the $(\mu_{i,j})_{i\leq n,j\leq k}$ are distinct. Thus we may replace, in Eq (\ref{commutat1}), the polynomial  $X^k+X^{k-1}B_{k-1}+\cdots+XB_1+B_0$ with a sparse polynomial. Yet, necessarily, $B_0$ must appear (cf. Eq (\ref{binome}) in Section 2) ; else $x=0$ is a common root of the polynomials $(\theta(x,\lambda_i))_i$. 
  \end{rem} 
   \begin{cor}
    Consider the Riccati Eq (\ref{comp}) where the coefficients $A,B_1,B_2,C$ are generic commuting $n\times n$ matrices. Then there are $2^n$ solutions, in $X\in\mathcal{M}_n(\overline{K})$, and any solution commutes with the coefficients.
  \end{cor}
  \begin{proof} 
  As in the proof of Theorem \ref{aston}, we show that $A$ is invertible. Using Section 3, a solution $X$ of  Eq (\ref{comp}) can be write $X=A^{-1}Z-A^{-1}B_2$ where $Z$ is a solution of the equation $Z^2+DZ+E=0_n$ with $D=-B_2+AB_1A^{-1}=-B_2+B_1,E=-AB_1A^{-1}B_2+AC=-B_1B_2+AC$. Note that $D,E$ are generic commuting matrices. According to Theorem \ref{aston}, $Z$ is a polynomial in $E$. Moreover $A,A^{-1},B_1,B_2$ are polynomials in $C$. Then $Z$ is a polynomial in $C$ and $X$ too.  
  \end{proof}
	Using Maple, we can prove, when $n\leq 3$, that any solution of the commuting Eq (\ref{plex1}) commute with $B_1,B_2,C$. \\ 
  Moreover, consider, for $n=3$, the equation of degree $4$
  \begin{equation}  \label{inst} X^4+A_1XA_2XA_3X+X^3+X^2+XA_4XA_5+A_6XA_7+B=0_3  \end{equation} 
  where the $3\times 3$ matrices $(A_i)_i,B$ are generic commuting matrices. Using Maple, we specialize Eq (\ref{inst}) and we obtain almost always that the solutions in $X\in\mathcal{M}_3(\overline{K})$ commute with $(A_i)_i,B$.
  Other numerical experiments do seem to indicate that the following surprising result is true
  \begin{conj}   \label{surp} Let $g$ be a non-zero polynomial over $K$ in $k+1$ non-commuting indeterminates $(u_1,\cdots,u_{k+1})$ such that\\
  $i)$ $g(u_1,\cdots,u_k,0)=0$\\
   and for every $i\leq k$\\
  $ii)$ $u_i$ appears exactly one time in the expression $g(u_1,\cdots,u_{k+1})$.\\
  $iii)$ The degree of $g$ with respect to $u_i$ is $1$.\\
  We consider the equation 
  \begin{equation} \label{surprise} g(A_1,\cdots,A_k,X)+B=0_n\text{ in the unknown }X\in\mathcal{M}_n(\overline{K}) \end{equation}
   where $A_1,\cdots,A_k,B$ are $n\times n$ generic commuting matrices. Then any solution of Eq (\ref{surprise}) commutes with the $(A_i)_i$ and $B$.  
  \end{conj}  
 The following shows that Conjecture \ref{surp} does not work if condition $ii)$ or $iii)$ is not fulfilled. 
  \begin{prop}
We consider the commuting Eq (\ref{plex2}): $X^2+BXB+C=0_2$, where $B,C$ are $2\times 2$ generic commuting matrices. Let $L=K(\sigma(B),\sigma(C))$. There are four trivial isolated solutions that commute with $B,C$. Moreover there are infinity many solutions $X\in\mathcal{M}_2(L(u))$ that depend on a parameter $u$ ; these solutions do not commute with $B,C$.  
   \end{prop}
  \begin{proof}
   We may assume that $B,C$ are diagonal matrices. Using Maple, we construct a Gr\"obner basis of the ideal generated by the associated $4$ algebraic equations in $4$ unknowns and $4$ indeterminates.  Since the Hilbert dimension of the ideal is $5$, the non-trivial solutions depend on $5-4=1$ parameter. 
    \end{proof}
    \begin{rem}
    More generally, the previous equation admits, in $\mathcal{M}_n(\overline{K})$, non-trivial solutions that depend on $n-1$ parameters.
    \end{rem}
 
  \section{Quadratic equation in matrices}
   
  The general quadratic matrix equation, in the unknown $X=[x_{i,j}]\in\mathcal{M}_n(\overline{K})$, has the form (cf. \cite{11})
  \begin{equation}   \label{degree2}
  \sum_{i=1}^r (A_iXB_iXC_i+D_iXE_i)+F=0_n  \end{equation}
  where $(A_i),(B_i),(C_i),(D_i),(E_i),F$ are $n\times n$ generic matrices.\\
    \indent Assume $n=2$. Eq (\ref{degree2}) is equivalent to $4$ algebraic equations of degree $2$ in the $4$ unknowns $(x_{i,j})_{i,j}$. Using a supplementary unknown $T$, we homogenize the previous equations and we seek the intersections of $4$ quadric hypersurfaces $(H_i)_i$ in the projective space $P_4(\overline{K})$. Let $\nu$ be the number of solutions of Eq (\ref{degree2}) in $P_4(\overline{K})$, counted with multiplicity, and by considering points at infinity. According to B\'ezout's Theorem,\\
  $i)$ either $\nu$ is finite and $\nu=2^4=16$. Moreover, if we know $16$ isolated solutions, then, necessarily $\nu$ is finite and there are no other solutions.\\
  $ii)$ or $\nu$ is infinite and there is an algebraic subset of solutions of positive dimension.\\
  \indent Note that the solutions such that $T=0$ satisfy the homogeneous equation
  \begin{equation}   \label{infty}
  \sum_{i=1}^r (A_iXB_iXC_i)=0_2.
  \end{equation}
   Thus Eq (\ref{degree2}) admits no solutions at infinity when Eq (\ref{infty}) admits the unique solution $X=0_2$. 
   \begin{exam2}  \label{infi} 
 $i)$ The set of solutions of the Riccati Eq (\ref{comp}) consists of $6$ isolated solutions in $\overline{K}^4$ and, at infinity, the set $\{X\not= 0\;|\; XAX=0_2\}=\{X\not=0\;|\;(XA)^2=0_2\}$ ; the previous set is a blunted cone of dimension $2$, that is a curve in the hyperplane $T=0$ of $P_4(\overline{K})$. \\
 $ii)$  Consider the quadratic equation $AXB_1X+XB_2X+DX+F=0_2$ where $A,B_1,B_2,D,F$ are $2\times 2$ generic matrices. Numerical experiments seem to indicate that this equation admits $8$ distinct solutions in $\overline{K}^4$. At infinity, we study the equation $\psi(X)=AXB_1X+XB_2X=0_2$. We obtain almost always $4$ distinct solutions in $P_4(\overline{K})$ ; if $X_0$ is such a solution, then $\dfrac{\partial \psi}{\partial X}(X_0)(X_0)=2\psi(X_0)=0_2$. Thus $rank(\dfrac{\partial \psi}{\partial X}(X_0))<4$ and $X_0$ is a multiple solution ; in order to obtain $16$ solutions in $P_4(\overline{K})$, its multiplicity would be $2$.  
  \end{exam2} 
   We prove \textbf{Theorem \ref{maxsol}} in which we give an instance of a generic quadratic matrix equation that admits $16$ isolated solutions in $\overline{K}^4$.
     \begin{proof}
	 $i)$  Let $\phi:X\rightarrow AXB_1X+XB_2X+X^2C+DX$. According to Theorem \ref{isol}, Eq (\ref{degmax}) has $\delta=\mathrm{degree}(\phi)$ solutions in $\mathcal{M}_2(\overline{K})$. According to B\'ezout's theorem, $\delta\leq 16$ and moreover, for every specialization of $F$, if the number $\nu$ of solutions of Eq (\ref{degmax}) is finite, then $\nu\leq \delta$ (cf. \cite[Section 6.3]{22}). \\
	 We randomly choose a specialization of Eq (\ref{degmax}) over $\mathbb{Q}$ and we study it with Maple. In general such an equation (just get one) admits a Gr\"obner basis in the form 
	 $$p_0(x_{1,1})=0\text{ where }p_0\in\mathbb{Q}[x]\text{ has degree }16\text{ and } S_{16}\text{ as Galois group over }\mathbb{Q}$$
	$$\text{ and for every }(r,s)\not=(1,1),\;x_{r,s}=p_{r,s}(x_{1,1})\text{ where }p_{r,s}\in\mathbb{Q}[x]\text{ has degree }15.$$ 
	Thus $\nu=16$ and we deduce that $\delta=16$.\\
	Let $p\in K[x]$ be the monic polynomial, the roots of which, are the solutions in $x_{1,1}$ with multiplicity. Since $p_0$ is irreducible, $p_0$ is a specialization of $p$ in $\mathbb{Q}$, $\mathrm{degree}(p)=16$ and $p$ is irreducible over $K$. If $G$ is the Galois group of $p$ over $K$, then $S_{16}$ is a subgroup of $G$ and $G=S_{16}$. Finally, $p$ has exactly $16$ simple roots and, to each root $\tilde{x}_{1,1}$, is associated exactly one solution $(\tilde{x}_{i,j})_{i,j}$ of  Eq (\ref{degmax}). That implies that the $(\tilde{x}_{i,j})_{(i,j)\not=(1,1)}$ are in $K[\tilde{x}_{1,1}]$ and solving  Eq (\ref{degmax}) is equivalent to solve $p=0$.\\
	$ii)$ Finally Eq (\ref{degmax}) has no solutions at infinity, that is, the homogeneous equation Eq (\ref{homog}): $AXB_1X+XB_2X+X^2C=0_2$ has the unique solution $X=0_2$ and consequently, this solution has multiplicity $16$. 
     \end{proof}
        Assume $n=3$. We look for a generic quadratic equation that is maximal in the sense that it has $2^9=512$ isolated solutions in $\mathcal{M}_3(\overline{K})$. Unfortunately Eq (\ref{degmax}) is not convenient. The set of non-zero solutions of Eq (\ref{homog}) is a cone in $\overline{K}^9\setminus\{0\}$. Therefore we consider solutions in $P_8(\overline{K})$. When we specialize Eq (\ref{homog}), the Hilbert dimension of the associated ideal is  almost always $1$, that is, there is at least one curve of solutions in $P_8(\overline{K})$. To obtain a maximal equation, we proceed as follows : first,  we seek a homogeneous equation $\psi(X)=0_3$ that has no non-zero solutions ; then, we consider the equation $\psi(X)+DX+F=0_3$ where $D,F $ are $3\times 3$ generic matrices. According to Theorem \ref{isol}, the previous equation has only isolated simple solutions ; thus it has $512$ distinct solutions in $\mathcal{M}_3(\overline{K})$.
  According to numerical experiments, we conjecture the following
      \begin{conj}
  The $3\times 3$ matrices $(A_i),(B_j),(C_k),D,F$ are generic. Then the equation in the unknown $X\in\mathcal{M}_3(\overline{K})$
  \begin{equation*}
  \psi(X)=A_1XB_1X+A_2XB_2X+XB_3XC_1+XB_4XC_2+XB_5X=0_3
  \end{equation*} 
  has no non-zero solutions, or equivalently, the solution $X=0_3$ has multiplicity $512$.\\
  Moreover the SC of the equation $\psi(X)+DX+F=0_3$ is $(512,S_{512})$.
  \end{conj}
  
  \section{Instances of Riccati equations}
  In this section, we adopt another point of view. Let $n=2$ and $\nu$ denote the number of solutions, \emph{with multiplicity} and not at infinity, of Eq (\ref{comp}). Here the matrices $A,B_1,B_2,C$ are numeric and fixed in $\mathcal{M}_2(\mathbb{C})$ ; the only condition that is required on these matrices is $A\not=0_2$. We look at the different forms of the set of solutions of a non-generic Riccati equation. The results that follow are known by V.V. Palin but, in the literature, they are partially treated, as in \cite{21} ; as a supplementary result, we show that, when there is an infinity of solutions, the algebraic sets of solutions have dimension $1$ or $2$.\\ 
  $\bullet$  We know that if $\nu$ is finite, then $\nu\leq 6$. For every $r\leq 6$, we give an instance of Eq (\ref{comp}) such that $\nu=r$. Another supplementary result is that, often, these instances, in the form $\phi(X)=Y$, are such that $Y$ is a ramification point of $\phi$. Recall that, if $M$ is non-derogatory, then $\nu$ is finite.\\
$i)$  $r=0$. The function $\phi_0(X)=X^2$ has degree $4$; an equation $\phi_0(X)=Y$ has $0,4$ or an infinity of solutions. In particular, the equation $\phi_0(X)=\begin{pmatrix}0&1\\0&0\end{pmatrix}$ has no solutions. Then $\phi_0$ has no ramification points.\\
$ii)$  $r=1$. The function $\phi_1(X)=X\begin{pmatrix}0&1\\0&1\end{pmatrix}X+\begin{pmatrix}-1&2\\0&-1\end{pmatrix}X$ has degree $2$. The equation  $\phi_1 (X)=Y$, where $Y\in\mathcal{R}=\{\begin{pmatrix}u&v\\-w&w\end{pmatrix}\;|\; w\not=0\text{ or }u+v\not=2\}$, has one solution because the unique convenient $M$-invariant subspace of dimension $2$ is $\ker(M^2)$. Thus the elements of $\mathcal{R}$ are ramification points of $\phi_1$.\\
  \emph{Question :} Does such a sole solution exist when $A$ is invertible ? \\
$iii)$  $r=2$. The function $\phi_2(X)=X^2+X$ has degree $4$, as every generic polynomial function of $\mathbb{C}[X]$. The equation $\phi_2(X)=Y=\begin{pmatrix}u&0\\-2&u\end{pmatrix}$ has $2$ solutions if $u\not=-1/4$. The matrices $Y$ are ramification points of $\phi_2$.\\
	$iv)$  $r=3$. The function $\phi_3(X)=X^2+\begin{pmatrix}2&0\\0&-1\end{pmatrix}X$ has degree $6$. The equation $\phi_3(X)=Y=\begin{pmatrix}-1&1\\0&2\end{pmatrix}$ has a unique solution $X_0$ that has multiplicity $3$ (see Example \ref{crit}, $i)$, in Section 4).
Here $\sigma(M)=\{1,1,1,-2\}$ and $M$ is non-derogatory. If we perturb slightly the matrices $Y$ (by adding small complex numbers), then, according to Theorem \ref{isol}, we obtain (in general) a generic equation that has $6$ isolated solutions in a neighborhood of $X_0$ and not $3$ solutions ; that is because $Y$ is a ramification point of $\phi_3$. \\
 $v)$ $r=4$. Use a generic commuting Eq (\ref{comp}) (cf. Theorem \ref{aston}). The equation $\phi_3(X)=Y$, where $Y$ is a generic diagonal matrix, has $4$ solutions. The generic diagonal matrices are ramification points of $\phi_3$.	Then $M$ has four distinct eigenvalues.\\
 $vi)$ $r=5$. The equation $\phi_3(X)=\begin{pmatrix}u&v\\0&w\end{pmatrix}$ has, for generic $u,v,w$, $5$ simple solutions. The generic matrices  $\begin{pmatrix}u&v\\0&w\end{pmatrix}$ are also ramification points of $\phi_3$.\\
$vii)$  $r=6$. Use a generic Eq (\ref{comp}) (cf. Theorem \ref{ric}) where $M$ has four distinct eigenvalues\\
  $or$ this instance, communicated to us by V.V. Palin. The function
	$$\phi_6(X)=X\begin{pmatrix}0&1\\1&1\end{pmatrix}X-\begin{pmatrix}0&1\\0&0\end{pmatrix}X+X\begin{pmatrix}0&1\\0&0\end{pmatrix}$$
	has degree $6$. The equation $\phi_6(X)=0_2$ admits the sole solution $X=0_2$, with multiplicity $6$.
Then $M$ is nilpotent and non-derogatory; remark that $\ker(M^2)$ is the unique $2$-dimensional $M$-invariant subspace. The previous example can be generalized as follows
\begin{rem}
Consider the equation $XAX+BX-XB=0_2$ where $A,B$ are $2\times 2$ generic matrices. Then we can show that the previous equation admits the trivial solution $X=0_2$ with multiplicity $4$ and $2$ supplementary non-zero solutions.
\end{rem}
 $\bullet$ Assume that $\nu$ is infinite. As we shall see, the Hilbert dimension $d$ of the ideal generated by the solutions is $1$ or $2$. 
 \begin{rem} In Section 3, we saw that, if $A$ is invertible, then Eq (\ref{comp}) can be rewritten in the form of Eq (\ref{final}) : $X^2+BX+C=0_2$. One has $\det(M-\lambda I_2)=\det(\lambda^2 I_2-\lambda B+C)$. Assume that $B,C$ commute ; then we may assume that $B,C$ are upper-triangular and there are orderings $(\mu_1,\mu_2)$ and $(\nu_1,\nu_2)$ of $\sigma(B)$ and $\sigma(C)$ such that the eigenvalues of $M$ are the roots of the polynomial $(\lambda^2-\mu_1 \lambda+\nu_1)(\lambda^2-\mu_2 \lambda+\nu_2)$. Since $\nu$ is infinite, the previous polynomial has a double root.
 \end{rem}
$i)$  $d=1$. Use a generic Eq (\ref{binome}) $X^2+TX=0_2$ (cf. Proposition \ref{grass}). \\
	Here $\sigma(M)=\sigma(T)\cup\{0,0\}$ ($3$ distinct eigenvalues) and $M$ is derogatory and diagonalizable. Moreover $\mathcal{S}$ (cf. Section 3) consists of $2$ non-intersecting straight lines in the direction of $2$ nilpotent matrices.\\
  \emph{Question :} Does there exist an instance of Eq (\ref{comp}) such that $d=1$ and that the variety of solutions is not flat ?\\
 $ii)$   $d=2$. Consider the equation $\phi_0(X)=I_2$. The non-singular solutions are $\pm I_2$ and the set of singular solutions is $\{\begin{pmatrix}a&b\\c&-a\end{pmatrix}\;|\;a^2+bc=1\}$, that is $2$ isolated points and a hyperboloid of one sheet of dimension $2$. Here $\sigma(M)=\{\lambda,\lambda,\mu,\mu\}$ where $\lambda\not=\mu$ and $M$ is diagonalizable and derogatory.\\
  $or$ $XAX=0_2$ where $A\not=0_2$.	Here $M$ is nilpotent and derogatory and the set of solutions is a cone of dimension $2$ when $A$ is invertible (see Example \ref{infi}. $i)$) and else, considering the cases $A=\begin{pmatrix}0&1\\0&0\end{pmatrix}$ and $A=\begin{pmatrix}0&0\\0&1\end{pmatrix}$, the union of two intersecting planes.\\
  The previous instance shows that the set of solutions at infinity is always a curve in $P_4$, that implies that the maximal Hilbert dimension of the ideal generated by the solutions, not at infinity, is $2$. We deduce the following
  \begin{prop}
  We consider Riccati Eq (\ref{comp}) where $n=2$ and $A\in\mathcal{M}_2(\mathbb{C})\setminus \{0_2\},B_1,B_2,C\in\mathcal{M}_2(\mathbb{C})$ are given fixed numeric matrices.
  There exist instances of the previous matrices such that the set of solutions, not at infinity, of the associated Eq (\ref{comp}) is one of the following \\
  $i)$ $0,1,\cdots$ or $6$ elements.\\
  $ii)$ An algebraic set of dimension $1$.\\
  $iii)$ An algebraic set of dimension $2$.\\
  Conversely, each instance of Eq (\ref{comp}) has a set of solutions that has one of the previous three forms.  
  \end{prop}

      \textbf{Acknowledgments}\\
The author thanks Michel Laurent and Daniel Barlet for many valuable discussions; he thanks the referee for his careful reading of the paper.

\bibliographystyle{plain}

\end{document}